\documentclass[10pt]{article}
\usepackage{amsmath,amssymb,amsthm}

\title{Lower bounds on directional complexity for irrational triangle billiards}
\author{Dmitri Scheglov\\
University of Oklahoma}

\theoremstyle{plain}

\begin{document}

\maketitle
\setlength{\parindent}{0pt}

\begin{abstract}
\noindent
We provide explicit lower bounds on directional complexity for a class of irrational triangle billiards  for a full measure $F_\sigma$-set of directions.
\end{abstract}

\section{Introduction and main results}

This paper is a continuation of our paper [7], where we provided an explicit version of the theorem by Galperin, Kruger and Troubeczkoy on the splitting of a thin parallel beam of triangle billiard trajectories. In paper [7] we gave an explicit upper bound on the splitting time in terms of a particular number theoretic function of angles $F_{\alpha\beta}(\epsilon)$. 
\

We also described some class of angles $\alpha, \beta$ for which  the function $F_{\alpha\beta}(\epsilon)$ is correctly defined and we conjectured that  $F_{\alpha\beta}(\epsilon)$ is in fact defined for any pair of irrational numbers $\alpha,\beta$ and moreover is uniformly bounded by some universal function $F(\epsilon)$, namely $F_{\alpha\beta}(\epsilon)\leq F(\epsilon)$.
\

\

Here we continue to exploit this function $F_{\alpha\beta}(\epsilon)$ and in particular we provide an explicit lower bound on the important dynamical characteristic of a billiard - its directional complexity.
\

\

The complexity of a polygonal billiard is a dynamical characteristic which roughly speaking measures the growth of combinatorial types of different orbits. Namely one enumerates the sides of the $k$-gon by symbols $1, 2,\ldots, k$ and then associates a word from this alphabet to a billiard trajectory of length $n$ reading the sides of the polygon which it hits.
\

The complexity function $p(n)$ is a total number of different words of length $n$ obtained this way.

\

There is also an analogous definition of the directional complexity function. We introduce a restriction on the trajectories, namely we consider the trajectories starting from a fixed side of a polygon under a particular angle $\theta$. And then we again count the number of different words of the length $n$. The resulting function is called a directional complexity in the direction $\theta$. We will denote it as $p_{\theta}(n)$ or sometimes just $p(n)$ if it does not lead to ambiguity.

\

\

By trivial reasons complexity function can not grow faster than exponentially and A. Katok [4] proved that for any polygon the complexity function in fact grows subexponentially, but his estimate is not explicit. It is still the best known upper bound on the complexity growth and it is a difficult open problem to provide any explicit subexponential upper estimate on $p(n)$.
\

\

S.Troubeczkoy found a quadratic lower bound for $p(n)$ in case of any polygon.[8]
\

\

In case of a rational polygon, meaning that the angles are rational multiples of $\pi$, the billiard is essentially equivalent to the geodesic flow on the flat compact surface with a finite number of singularities. Then one can use Teichmuller theory to investigate the billiard dynamics.
\

In this case using the results of H.Masur on the growth for quadratic differentials [6], J.Cassaigne, P. Hubert, S. Troubeczkoy [1] proved that there are constants $c_{1}$ and $c_{2}$ such that $c_{1}n^{3}<p(n)<c_{2}n^{3}$.
\

For a directional complexity $p_{\theta}(n)$ of any polygon with $k$ sides there is a polynomial upper bound by  E.Gutkin and S. Troubeczkoy [3]:
\

$p_{\theta}(n)\leq krn{\bigl (1+\frac{n}{2}\bigr )}^{s}$,

\

where $r$ is the least common denominator of the rational angles of $P$ and $s$ is a number of distinct irrational angles of $p$.
\

Note that this upper boundary is quite universal, namely it is basically independent on the angles of the polygon and on the direction.
\

The natural question is then to give the $lower$ bound on the directional complexity.  This question however is more delicate and the current paper is devoted to the lower bound estimates on directional complexity for irrational triangles.
\

\section{Intuitive scheme of the proof.}
\

First of all let us consider a simple example of a biliard in the square and horizontal direction $\theta$. By " horizontal " we mean the direction parallel to one of the square sides. In this direction $p(n)=1$ which means there is no growth at all. The same picture happens in any rational direction in the square, for which $p(n)$ is bounded. 
\

Intuitively this happens because any orbit in a rational direction is periodic and so the whole billiard flow in such a direction just splits into a finite number of "periodic beams" without producing any complexity.
\

On the other side the complexity in any irrational direction is closely connected to the so-called " Sturmian sequences" and certainly grows.
\

\

In order to understand the nature of the complexity growth we make the folowing well-known observation. Namely, the function $p_{\theta}(n)$ increases at those time moments $n_{0}$ when the thin beam of parallel trajectories of a given combinatorics, starting in a given direction, hits a polygon vertex.
\

 At this time moment one " half" of the beam hits one side of the polygon and another " half" of the beam hits another side and so we have at least two different words of length $n_{0}+1$ which implies $p_{\theta}(n_{0}+1)\geq p_{\theta}(n_{0}) +1$.
\

\

This mechanism of complexity growth however fails in general if the direction $\theta$ is periodic. In this case there exist thin " periodic beams" which never split and do not produce complexity. The before mentioned example of a rational direction in the square billiard perfectly demonstrates this fenomenon.
\

It is theoretically possible that for some irrational polygons there exist periodic directions where not all the orbits are periodic. In the other words so-called " partially - periodic directions" may exist. However the analysis of periodic trajectories in irrational polygons so far is very difficult problem and it is not even known if periodic orbits always exist. 
\

On the other side even if almost all the orbits in a given direction $\theta$ are periodic, it could possibly happen that there is an infinite number of periodic beams, which would imply the growth of $p_{\theta}(n)$. But again, for generic polygons there are yet no tools, allowing one to effectively analyse the behaviour and distribution of periodic orbits.
\

\

Having this said it is natural to consider a direction $\theta$ which is not periodic, in the other words, there are no periodic orbits in $\theta$-direction. A priory there is still a possibility that the thin parallel beam of trajectories never splits, however this possibility is excluded by the very nice result of Galperin, Kruger and Troubeczkoy[2], which can be reformulated as follows:
\

\

\textbf{Theorem ( Galperin, Kruger, Troubeczkoy).} In any polygonal billiard any non-parallel beam of trajectories splits.
\

\

We would like to remark that the theorem above in the original paper was formulated slightly differently but the given formulation easily follows.
\

Having in mind this result the approximate scheme on the estimating $p_{n}(\theta)$ would approximately look as follows.
\

We divide a side of the polygon on $m$ equal pieces of length $1/m$ and start a parallel beam from each piece. If there is a uniform splitting time time $T(m)$ for each beam, then obviously $p(T(m))\geq m $. The numbers $T(m)$ form a very special increasing sequence, so in order to estimate $p(n)$ we find a maximal $m$ such that $T(m)\leq n$, which provides a lower bound $p(n)>m$.

\

However we would like to explain the reader several important obstacles which we must overcome during this process.
\

\

First of all even if the beam $B(\epsilon)$ is not periodic, it can be very close to a periodic beam. Namely if $\theta$ is a beam direction and $\alpha$ is a direction of a fixed periodic beam, than if $\epsilon$ is small enough and the middle point of the base of the beam $B(\epsilon)$ coincides with a middle point of periodic beam, then as the value $\phi=|\theta-\alpha|\longrightarrow 0$ implies that the splitting time $n_{0}(\epsilon)\longrightarrow\infty$.
\

\

Informally speaking this example shows that the periodic orbits are " the obstacles" to the uniform splitting.
\

\

And the second issue to keep in mind is the following. Assume that our beam $B_{\theta}(\epsilon)$ is in the reasonable sense "far" from periodic orbits. The precise meaning of the word "far" will be explained little later. We then have to provide an explicit bound for the splitting time $n_{0}(\epsilon)$ in the other words an effective version of the splitting theorem by Galperin, Kruger and Troubeczkoy.
\

\

The main result of the paper [7] was a dichotomy, briefly formulated as follows: either beam $B(\epsilon)$ contains a periodic orbit inside or its length is uniformly bounded by some function $M(\epsilon)$. In the other words " parallel $\epsilon$ - beams, not containing periodic orbits can not be too long".
\

\

Having this dichotomy we may assume that if the beam $B(\epsilon)$ does not split in time $M(\epsilon)$ then it has a periodic orbit inside. For any fixed direction $\theta$ we introduce a sequence of functions $\phi (N)=min |\theta-\alpha|$, where $\alpha$ runs through all periodic directions, corresponding to the directions of the length less than $N$.
\

Using some delicate analysis we estimate the splitting time for a beam which possibly contains a periodic orbit inside, in terms of the sequence $\phi_{n}$.
\

We would like to note that this result is independent on the fact whether or not there exist periodic orbits inside a triangle. If in some triangle there are no periodic orbits at all ( which presumably never happens), then the corresponding terms in the explicit formulas for the splitting time would vanish.
\

\

This way we obtain the explicit formulas for any non-periodic directions $\theta$ in terms of the sequence $\phi_{n}$. However these lower bounds are in a sense " not explicit enough", as we basically have no  information about the sequence $\phi_{n}$.
\

\

However these explicit formulas combined with some simple measure-theoretic arguments  allow us to provide a lower bound on the $p_{\theta}(n)$ for $typical$ directions $\theta$ as we are able to say that typical direction " stays away" from periodic directions with a prescribed distance.
\

\

Now it is time to turn into formal definitions.
\

\section{Definitions.}
\

In this section we borrow several definitions from the paper [7] in order to keep the current paper self-contained and provide the results needed for the proof of the main theorem.
\

\

First of all we will remind the so-called Katok-Zemlyakov construction[5], which is a main technical tool in our analysis of the billiard trajectories. Consider a billiard trajectory inside a polygon and in particular a moment when it hits the polygon side. Instead of reflecting the trajectory, we reflect a polygon, using the " optic law", that the angle of reflection coincides with the angle of incidence. 
\

From the point of view of the observer inside the polygon, looking at the side as in the mirror, the billiard trajectory instead of reflection goes further " behind the mirror". Then we continue the process indefinitely many times.
\

As a result,  instead of the piecewise linear billiard trajectory in a fixed polygon, we have a fixed straight line and a sequence of reflected polygons along the line.
\
In the picture below we show the result of the one-time application of the Katok-Zemlyakov construction to the triangle. This way we obtain a nice shape, which we call a $Kite$.

\begin{minipage}[c]{10mm}
\unitlength=1.00mm
\special{em:linewidth 0.4pt}
\linethickness{0.4pt}
\begin{picture}(23.67,44.67)
\put(13.00,44.67){\line(0,-1){41.67}}
\put(13.00,3.00){\line(-2,5){10.67}}
\put(2.34,30.00){\line(3,4){10.67}}
\put(13.00,44.34){\line(3,-4){10.67}}
\put(23.67,30.34){\line(-2,-5){11.00}}
\end{picture}

\end{minipage}
\

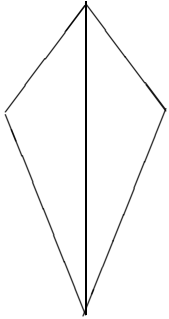
\textbf{Pic. 1.} Kite
\

\

The next picture represents Katok-Zemlyakov construction, applied several times to a kite.
\

\begin{minipage}[c]{10mm}
\unitlength=1.00mm
\special{em:linewidth 0.4pt}
\linethickness{0.4pt}
\begin{picture}(94.67,42.67)
\put(9.00,3.67){\line(-1,4){6.67}}
\put(2.33,29.67){\line(2,3){6.33}}
\put(8.66,39.34){\line(2,-3){6.67}}
\put(15.33,29.34){\line(-1,-4){6.00}}
\put(9.00,3.67){\line(1,1){17.33}}
\put(26.33,21.00){\line(0,1){13.00}}
\put(26.33,34.00){\line(-5,-2){11.33}}
\put(4.00,23.34){\vector(1,0){90.67}}
\put(26.33,34.00){\line(6,-1){11.00}}
\put(37.33,32.34){\line(1,-3){8.00}}
\put(45.33,8.34){\line(1,5){5.33}}
\put(45.33,8.34){\line(-3,2){19.00}}
\put(50.66,35.00){\line(-1,1){7.67}}
\put(43.33,42.67){\line(-3,-5){6.33}}
\end{picture}

\end{minipage}

\

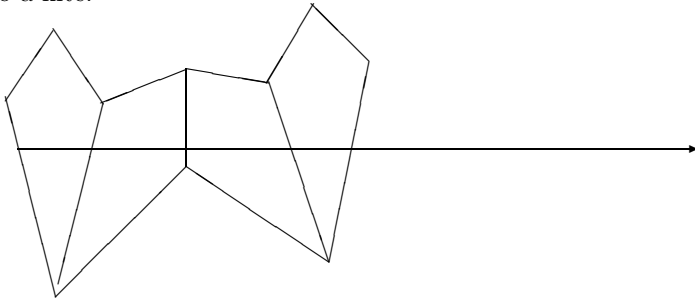
\textbf{Pic.2.} Katok-Zemlyakov construction applied to a kite.
\

\

We will also need the notion of the parallel beam of the trajectories and we borrow corresponding definition from [1].
\

\

\textbf{Definition 1.} An $(\epsilon, T)$ - beam is a set of parallel segments, corresponding to the application of the Katok-Zemlyakov construction along some direction and from some base point, where $\epsilon$ is a width of the beam and $T$ is the length of the maximal parallel segment. 
\

The left interval on the kite side, transversal to the beam direction is called \underline {a base segment} or  \underline {base of a beam} and the right interval is called \underline{ an end segment} or \underline{end of the beam}.
\

\

Note, that by the definition beam does not have any kite vertices inside, as the Katok-Zemlyakov is undefined on when the trajectory hits a vertex.
\

We will usually denote $(\epsilon, T)$ - beam as $B(\epsilon, T)$.
\

Below one may see a picture, which provides a geometric intuition behind the notion of $(\epsilon, T)$ - beam.
 \

\begin{minipage}[c]{10mm}
\unitlength=1.00mm
\special{em:linewidth 0.4pt}
\linethickness{0.4pt}
\begin{picture}(108.67,44.00)
\put(9.00,3.67){\line(-1,4){6.67}}
\put(2.33,29.67){\line(2,3){6.33}}
\put(8.66,39.34){\line(2,-3){6.67}}
\put(15.33,29.34){\line(-1,-4){6.00}}
\put(9.00,3.67){\line(1,1){17.33}}
\put(26.33,21.00){\line(0,1){13.00}}
\put(26.33,34.00){\line(-5,-2){11.33}}
\put(4.00,23.34){\vector(1,0){90.67}}
\put(26.33,34.00){\line(6,-1){11.00}}
\put(37.33,32.34){\line(1,-3){8.00}}
\put(45.33,8.34){\line(1,5){5.33}}
\put(45.33,8.34){\line(-3,2){19.00}}
\put(50.66,35.00){\line(-1,1){7.67}}
\put(43.33,42.67){\line(-3,-5){6.33}}
\put(91.67,12.00){\line(1,4){6.67}}
\put(3.00,27.00){\vector(1,0){92.67}}
\put(91.67,12.00){\line(4,5){17.00}}
\put(106.67,44.00){\line(-5,-3){8.67}}
\put(106.67,44.00){\line(1,-6){2.00}}
\put(108.67,32.67){\line(0,0){0.00}}
\end{picture}

\end{minipage}

\

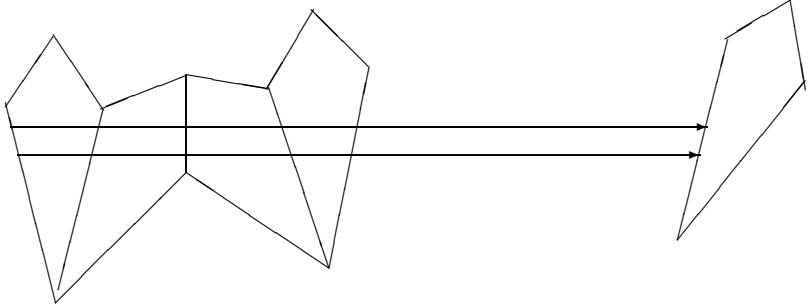
\textbf{Pic.3.} $(\epsilon, T)$ - beam of parallel trajectories.

\

We also remind the following well-known inequality, which will be used later.
\

\

\textbf{Lemma 1.} Let $S$ be a billiard trajectory, $L(S)$ be its geometric length and $N(S)$ be its combinatorial length, meaning the number of reflections in the Katok-Zemlyakov construction. Then there exist positive constants, $c$, $C$, depending only on the kite, such that: $cN(S)\leq L(S)\leq CN(S)$.
\

\

The proof of Lemma 1. is elementary and we will not reproduce it here.
\

From here and further we will denote as $C$ any large enough or small enough constant, depending on the context.
\

We will also use the following simplifying convention from [1], which would allow us to avoid  too complicated expressions:
\

\textbf{Convention.} Any function defined on the positive integers $f: \mathbb{Z}_{+}\rightarrow\mathbb{R}$ is by default extended to a function $f:\mathbb{R}_{+}\rightarrow\mathbb{R}$ by the rule $f(x)=f(\left [ x\right ] +1)$, where $x\in\mathbb{R}_{+}\setminus\mathbb{Z}_{+}$.

\

\
The next several definitions are borrowed from [1] in order to formulate a dichotomy theorem. We introduce them in order to keep our exposition self - contained.
\

\

\textbf{Definition 1.} Consider a finite subset $S\subset\Delta$ of a segment $\Delta$ of the circle. $S$ is called a relative $\epsilon$-net if after linear "blowing up" of the segment $\Delta$ to the length 1, $S$ becomes an $\epsilon$-net in the standard sence.
\

\

\textbf{Definition 2.} Let $0<\alpha, \beta<\pi$. A finite sequence $x_{1},\ldots, x_{n}\in\mathbb{S}^{1}$ is called $\alpha\beta$-connected if for any $i: 1\leq n-1$ we have that either $x_{i+1}-x_{i}=\pm\alpha$ or $x_{i+1}-x_{i}=\pm\beta$ .
\

\

For any finite set of points $S\subset\mathbb{S}^{1}$  let $|S|$ denote its cardinality. We then have the following definition.
Fix a pair of numbers $\alpha$, $\beta$ as above. The $net$-function $ F_{\alpha \beta}$: (0,1)$\rightarrow Z_{+}$ is defined as follows.
\

\

\textbf{Definition 3.} $ F_{\alpha \beta}\left (\epsilon\right )=min \lbrace n\in\mathbb{Z}_{+}| $ $\forall  \alpha\beta $ -connected $S\subset\mathbb{S}^{1}$ $,$ $|S|>n$  $\exists $ $\overline S \subset S$, $\overline S$$-$ relative $\epsilon$-net  $\rbrace$
\

\

Informally speaking $ F_{\alpha \beta}\left (\epsilon\right )$ is a minimal cardinality of $\alpha\beta$ - connected sequence which guarantees that it contains a relative $\epsilon$-net.

\

\textbf{Definition 4.} Let $\alpha,\beta$ be rationally independent. Then $N_{\alpha\beta}\left ( k\right ) =$ min $\lbrace  \langle n\alpha+m\beta\rangle |$ for all $ |n|+|m|\leq k| \rbrace$, where $\langle x\rangle$ is a distance from $x$ to the closest integer.
\

\

The function $F_{\alpha\beta}(\epsilon)$ defined above turns out to be extremely important for our purposes. In fact our approach works precisely for all angles $\alpha, \beta$ for which $F_{\alpha\beta}$ is defined. In the paper [1] we proved the existence and gave an explicit upper estimate on the function $F_{\alpha\beta}$  for all angles $\alpha, \beta$ for which $\beta$ is a fixed irrational number and $\alpha/\beta$ allows very fast approximation by rational numbers. This result shows that the class of angles, for which $F_{\alpha\beta}$ is explicitely defined is quite non-trivial.
\

In the paper [7] we conjecture that $F_{\alpha\beta}(\epsilon)$ is defined for all $\alpha, \beta$ and moreover, that it is uniformly bounded from above for all $\alpha, \beta$, namely $F_{\alpha\beta}(\epsilon)\leq F(\epsilon)$.
\

\

The proof of this number-theoretic conjecture seems to be quite important for understanding the triangle billiard dynamics. In particular due to the result of the current paper it would automatically give the explicite lower bounds on the typical directional complexity for all triangle billiards.
\

\

Even if the number-theoretic conjecture is not true for all pairs $\alpha, \beta$ its proof for some particular pair of angles $\alpha, \beta$ implies the lower bound on the directional complexity for that particular triangle.
\

\

We now formulate a theorem from [7] which will serve as a key tool in our approach.
\

\

\textbf{Theorem (Effective dichotomy).}  Let $\alpha$, $\beta$ be a pair of rationally independent numbers such that $ F_{\alpha \beta}$($\epsilon$) is a correctly defined function for positive $\epsilon$. And let  $ K_{\alpha \beta}$ be the kite of diameter 1 with angles  $\alpha$, $\beta$. Let B($\epsilon$, $M$) be a parallel beam. 
\

Let us also introduce the following notations: $P_{\alpha\beta} (\epsilon) = \frac{16}{\epsilon} F_{\alpha\beta} \bigl ( {\bigl ( \frac{\epsilon}{1600} \bigr )}^{\frac{16}{\epsilon}}            \bigr )$
\

Then either $M\leq M(\epsilon)= \frac{C}{N_{\alpha\beta} (P_{\alpha\beta}(\epsilon))}$ or $B$ contains a periodic trajectory inside, starting from the base. Here $C$ is a constant, depending only on the kite.
\

\

One more definition precisely defines the speed of approximation of a given direction by periodic.
\

\

\textbf{Definition 7.} Fix a kite $K$ on the plane, choose a side of $K$ and a direction $\theta$. As a tangent bundle to $K$ is naturally trivialized, we may think of any direction as of the point on the circle $\mathbb{S}^{1}$. Then we define a sequence $\phi(n)$ as follows:
\

\

$\phi(n)= min\lbrace|\theta-\alpha|\rbrace$, where $\alpha$ runs through all the periodic directions of the combinatorial length less than $n$.
\

\

Here by periodic direction $\alpha$ we mean a direction such that there exists at least one periodic orbit starting from the given side in the direction $\alpha$.
\

We should notice that the set of directions $\alpha$ from the definition above is finite and so the minimum is taken over a finite number of values. It happens because any periodic direction is uniquely defined by the combinatorics of a corresponding periodic trajectory which  is an easy well-known observation.
\

\

In the hypothetical case when there are no periodic trajectories at all, we put $\phi(n)=1$ identically.

\section{Complexity growth.}
\

We  now have all the tools to provide the lower bound for a directional complexity growth. First we prove a useful theorem, which estimates the splitting time of a thin beam of width $\epsilon$ in the direction $\theta$.
\

\

\textbf{Theorem.} Let $\alpha, \beta$ be a pair of kite angles, such that the function $F_{\alpha\beta}(\epsilon)$ is correctly defined. Let us also fix a side of the kite $I$, a non-periodic direction $\theta$ and $\epsilon$ small enough. Let $\phi(n)$ be the sequence, defined above and $B_{\theta}(\epsilon, T)$ be a parallel beam in the direction $\theta$.
\

Note that non-periodicity of $\theta$ implies that all the terms $\phi(n)$ are positive. 
\

Then $T\leq T(\epsilon)= max\lbrace M(\epsilon), \frac{C}{\phi(C\cdot M(\epsilon)) }\rbrace$, $C>1 $ large enough.
\

\

\textbf{Proof.} Let us assume that $T>T(\epsilon)$, which implies $T>M(\epsilon)$.  By the theorem 1 there is a periodic trajectory $S$ of the length $L(S)$ inside the beam $B_{1}=B_{\theta}(\epsilon,M(\epsilon))$, starting from its base. It implies that there is a subbeam $B^1_\theta = B_{\theta}(\epsilon, L(\epsilon))$ with the same base, which has parallel base and end kites and contains a small parallel beam of periodic orbits close to $S$.
\

We are going to take a closer look on the trajectory $S$ and estimate how the endpoints of $S$ are located inside $B_{1}$.
\

\begin{minipage}[c]{10mm}
\unitlength=1.00mm
\special{em:linewidth 0.4pt}
\linethickness{0.4pt}
\begin{picture}(141.00,114.66)
\put(19.33,4.33){\line(-1,4){17.00}}
\put(19.33,4.66){\line(1,4){17.00}}
\put(36.33,72.33){\line(-1,1){17.00}}
\put(19.33,89.33){\line(-1,-1){17.00}}
\put(124.00,29.66){\line(-1,4){17.00}}
\put(124.00,30.00){\line(1,4){17.00}}
\put(141.00,97.66){\line(-1,1){17.00}}
\put(124.00,114.66){\line(-1,-1){17.00}}
\put(11.00,37.66){\line(1,0){111.00}}
\put(11.00,37.66){\line(4,1){104.67}}
\put(3.00,69.33){\line(1,0){111.33}}
\put(114.33,69.33){\line(-4,-1){104.67}}
\put(10.00,41.66){\line(4,1){104.67}}
\put(10.33,40.00){\line(4,1){104.67}}
\put(10.66,38.66){\line(4,1){104.67}}
\put(4.67,63.67){\vector(1,0){8.00}}
\put(6.33,56.33){\vector(1,0){8.33}}
\put(8.00,49.67){\vector(1,0){8.67}}
\put(9.33,43.67){\vector(1,0){9.33}}
\put(109.67,87.00){\vector(1,0){12.67}}
\put(111.33,79.67){\vector(1,0){13.33}}
\put(113.67,71.33){\vector(1,0){13.33}}
\put(107.33,97.33){\line(1,0){18.33}}
\put(115.67,63.67){\line(1,0){14.67}}
\put(59.33,29.67){\makebox(0,0)[cc]{$L(\epsilon)$}}
\put(4.33,49.00){\makebox(0,0)[cc]{$\epsilon$}}
\put(48.67,40.67){\makebox(0,0)[cc]{$\alpha$}}
\put(115.33,50.33){\makebox(0,0)[cc]{$\delta$}}
\put(82.67,51.67){\makebox(0,0)[cc]{$L(S)$}}
\end{picture}

\end{minipage}
\

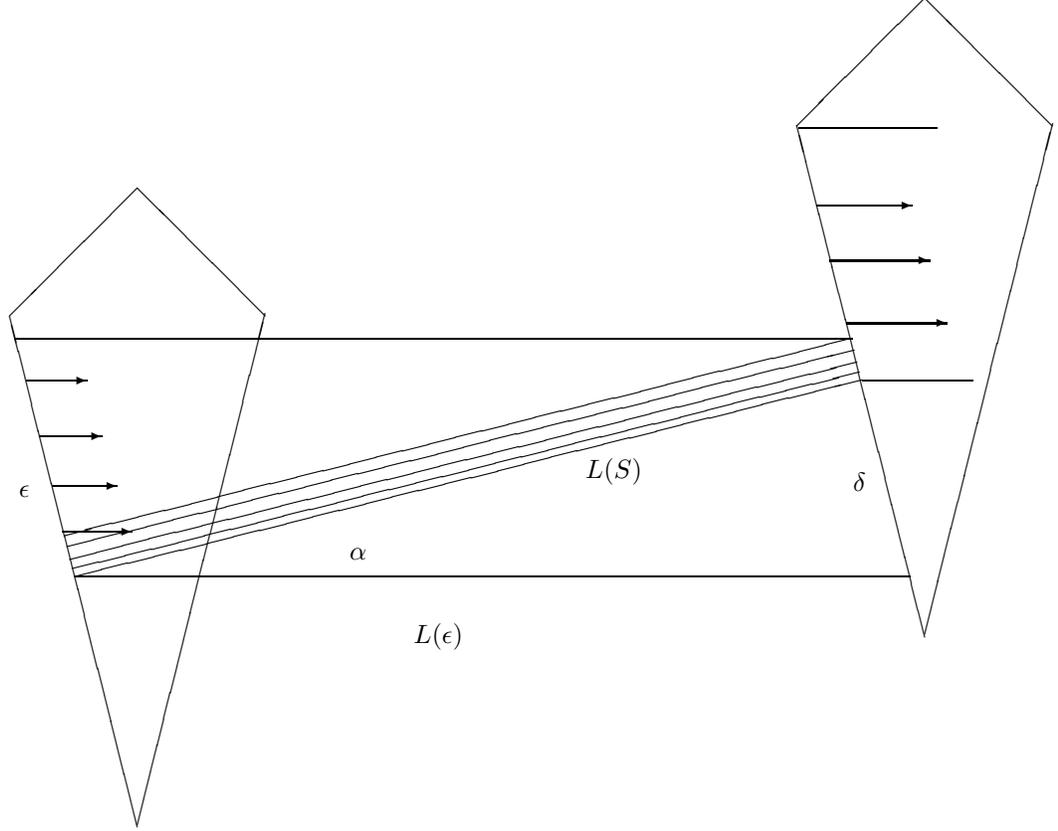
\textbf {Pic.4.} The beam $B^1_\theta = B_{\theta}(\epsilon, L(\epsilon))$ and a beam of periodic trajectories inside.

\

As we see on the picture there is a horizontal sheer on the length $\delta$ of the beam base. $\delta$ is exactly the difference between $\epsilon$ and the width of the periodic beam, located inside $B_{\theta}(\epsilon, L(\epsilon))$.
\

\

Let us now estimate the sheer $\delta$ from below.
\

As the beam length is $L(\epsilon)$ and $\alpha$ is exactly the angle between the periodic direction and $\theta$ then by definition of the sequence $\phi(n)$ we have $\alpha\geq \phi(CL(\epsilon))$. As for small enough $\alpha$ we have $ sin (\alpha )\approx\alpha$ then from the picture 4 we have $\delta\geq C sin(\alpha)L(\epsilon)\geq c\cdot sin(\phi (C\cdot L(\epsilon)))L(\epsilon)\geq c \phi(C\cdot L(\epsilon))L(\epsilon)\geq c\phi(C\cdot M(\epsilon))L(\epsilon)$.

\

\

We now take a new beam $B^{2}_{\theta}$ which base is a union of the base and end of the beam $B_{\theta}(\epsilon, L(\epsilon))$. The width of $B^{2}_{\theta}$ is $\epsilon + \delta$ which is clear from the picture 5. The picture 5 shows how the part of the end of the beam $B_{\theta}(\epsilon, L(\epsilon))$ is attached to its base.
\

\begin{minipage}[c]{10mm}
\unitlength=1.00mm
\special{em:linewidth 0.4pt}
\linethickness{0.4pt}
\begin{picture}(141.00,114.66)
\put(19.33,4.33){\line(-1,4){17.00}}
\put(19.33,4.66){\line(1,4){17.00}}
\put(36.33,72.33){\line(-1,1){17.00}}
\put(19.33,89.33){\line(-1,-1){17.00}}
\put(124.00,29.66){\line(-1,4){17.00}}
\put(124.00,30.00){\line(1,4){17.00}}
\put(141.00,97.66){\line(-1,1){17.00}}
\put(124.00,114.66){\line(-1,-1){17.00}}
\put(11.00,37.66){\line(1,0){111.00}}
\put(11.00,37.66){\line(4,1){104.67}}
\put(3.00,69.33){\line(1,0){111.33}}
\put(114.33,69.33){\line(-4,-1){104.67}}
\put(10.00,41.66){\line(4,1){104.67}}
\put(10.33,40.00){\line(4,1){104.67}}
\put(10.66,38.66){\line(4,1){104.67}}
\put(4.67,63.67){\vector(1,0){8.00}}
\put(6.33,56.33){\vector(1,0){8.33}}
\put(8.00,49.67){\vector(1,0){8.67}}
\put(9.33,43.67){\vector(1,0){9.33}}
\put(109.67,87.00){\vector(1,0){12.67}}
\put(111.33,79.67){\vector(1,0){13.33}}
\put(113.67,71.33){\vector(1,0){13.33}}
\put(107.33,97.33){\line(1,0){18.33}}
\put(115.67,63.67){\line(1,0){14.67}}
\put(59.33,29.67){\makebox(0,0)[cc]{$L(\epsilon)$}}
\put(4.33,49.00){\makebox(0,0)[cc]{$\epsilon$}}
\put(48.67,40.67){\makebox(0,0)[cc]{$\alpha$}}
\put(115.33,50.33){\makebox(0,0)[cc]{$\delta$}}
\put(17.33,12.33){\vector(1,0){36.33}}
\put(15.67,19.67){\vector(1,0){33.00}}
\put(14.00,26.33){\vector(1,0){30.67}}
\put(12.33,32.00){\vector(1,0){30.33}}
\put(10.67,22.67){\makebox(0,0)[cc]{$\delta$}}
\put(82.00,51.00){\makebox(0,0)[cc]{$L(S)$}}
\end{picture}

\end{minipage}

\

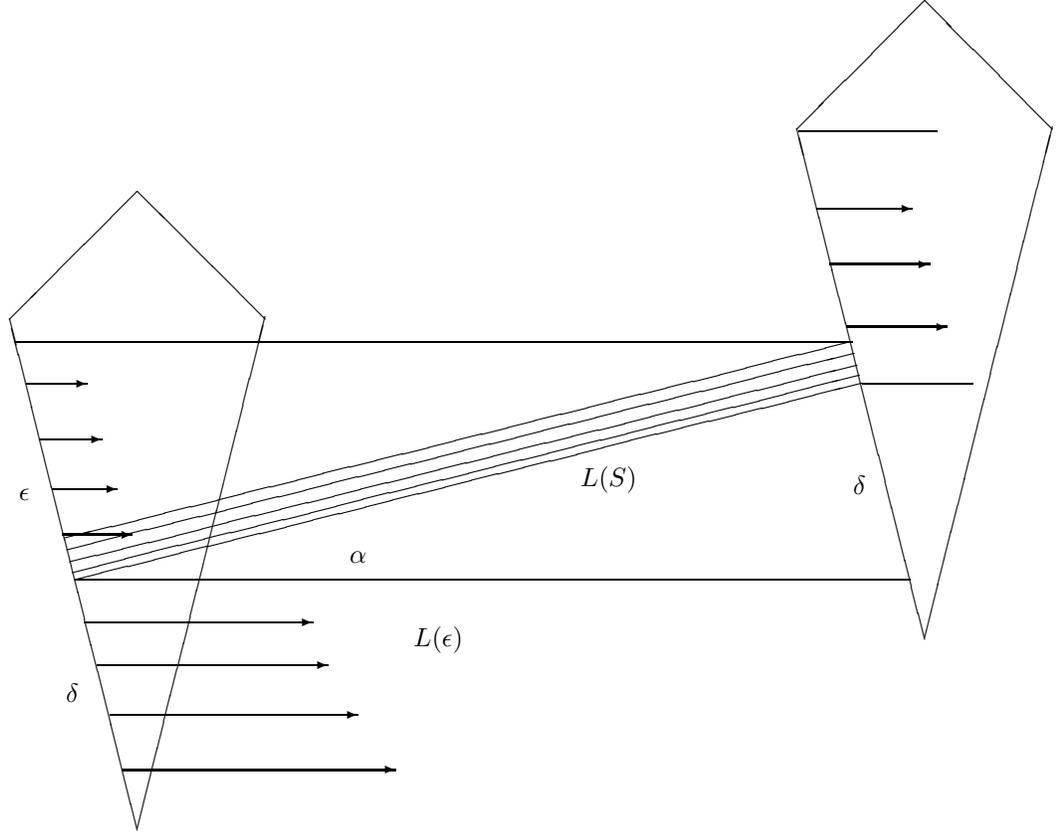
\textbf{Pic.5.} The second beam $B^2_\theta$ of the width $\epsilon+\delta$.

\

\

Assume first that the extended beam $B_{\theta}^{2}$ does not split at the time $L(\epsilon)$. 
\

In this case we still have a periodic trajectory $S$ inside the extended beam with angle $\alpha$ to the direction $\theta$. It implies that the previous argument can be applied to the extended beam $B_{\theta}(\epsilon, L(\epsilon))$ and we  construct a third beam $B^3_\theta$ and contiinue the process.
\

Note that each time the width of a new beam is greater then the width of the previous on  $\delta$.
\

\

Let us assume that this process terminates at $n$-th step, meaning that the $n$-th extended beam splits. As each of the extended beams has length less than $L(\epsilon)$ it implies that the splitting time $T(\epsilon)$ of the original beam satisfies $T(\epsilon)\leq  nL(\epsilon)$.
\

Since the perimeter of the kite is 1 and each time the beam width increases by $\delta$, it implies $\delta n\leq 1$.
\

Combining all the inequalities together we get: $T(\epsilon)\leq C/\phi(M(\epsilon))$.
\

\

As the computations above were maid under assumption that $L(\epsilon)\leq M(\epsilon)$ then in general for arbitrary beam $B_{\theta}(\epsilon, T)$ we have $T\leq  max(M(\epsilon), 1/\phi(C\cdot M(\epsilon)))$.
\

\

Let us now estimate the splitting time in the generic direction. First of all we remind an easy upper estimate on the number $P(n)$ of periodic beams of length less than $n$. As each such a beam is uniquely determined by its coding then the number of the periodic beams of length precisely $n$ is at most $4^{n}$. And sinse we count all the trajectories of the smaller length, we have:
\

$P(n)\leq 1+ 4 +\ldots +4^{n}\leq 4^{n+1}$
\

We now take a small enough $\rho>0$  and consider  an open set $U(n,\rho)= \lbrace \theta\in\mathbb{S}^{1}: |\theta-\alpha|< \rho 4^{-n}4^{-n-1}\rbrace$ where $\alpha$ runs through all periodic directions of length less then $n$.
\

As the number of periodic directions of length less then $n$ is bounded by $4^{n+1}$ we have that the Lebesguse measure of our set $ Leb(U(n,\rho))\leq\rho 4^{-n}$.
\

Let us denote $U(\rho)=\bigcup_n U(n,\rho)$

Summing up these inequalities for all $n$ we have $Leb(U(\rho))\leq\rho$.
\

\

Now for any direction $\theta$ in the complement of $U(\rho)$ we have by definition $\phi(n)\geq \rho 4^{-n}4^{-n-1}$ which implies that there is a large enough constant $C$ depending on $\rho$ that the splitting time of the parallel beam in the direction $\theta$ satisfies $T(\epsilon)\leq C \exp (C\cdot M(\epsilon))$ for large enough $C$
\
\

\

Each $U(\rho)$ is an open set so the complement $V(\rho)$ is closed. Let us now take a sequence $\rho_{n}=1/n$ and consider $V=\bigcup_n V(\rho_n)$.
\

The set $V$ belongs to the class $F_{\sigma}$ which is a countable union of closed sets and from considerations above easily follows:
\

$1)$ $V\in F_{\sigma}$ 
\

$2)$ $V$ has a full Lebesgue measure.
\

$3)$ For any direction $ \theta\in V$ we have $T(\epsilon)\leq C \exp (C\cdot M(\epsilon))$ for large enough $C$
\

\

Now let us turn to the directional complexity estimate. We first remind the idea in more details. 
\

Let us consider the billiard side $I$ and let all the the points of $I$ move in the direction $\theta$. Without loss of generality let us assume that the length of $I$ equals to 1.  We fix a positive integer $m$ on divide $I$ onto $m$ equal intervals of the length $1/m$.
\

Each of these intervals splits at the time $T(1/m)$, which in particular implies that $p_{\theta}(T(1/m))\geq m$.
\

Now we define the following function $N_{\theta}(n)=max\lbrace m | T(1/m)\leq n\rbrace$. As $p_{\theta}(n)$ is increasing function, the inequality above implies our main result, the lower bound on the directional complexity: $p_{\theta}(n)\geq N_{\theta} (n)$.
\

\

Let us finally summarize our observations to the following concluding theorem.
\

\

\

\

\textbf{Theorem 2. Lower bound on directional complexity.}
\

\

Let $\alpha, \beta$ be a pair of rationally independent irrational numbers, such that the function $F_{\alpha\beta}(\epsilon)$ is correctly defined. 
\

Let $R(m)= N_{\alpha\beta}\bigl ( 16m F_{\alpha\beta}\bigl ({\bigl (\frac{1}{1600m}\bigr )}^{16m}\bigr )\bigr )$ and let
\

\

$L(n)=max\lbrace m\in\mathbb{Z}_{+}| R(m)\geq C/\ln(n)\rbrace$, where $C$ is a big enough constant depending on $K$ and $\theta$.

Then for a kite $K$ with angles $\alpha, \beta$ and for any direction $\theta$ from $F_\sigma$-set of a full measure, the directional complexity $p_\theta (n)$ satisfies:
\

\

$p_\theta (n)\geq L (n)$

\end{document}